\documentclass[a4paper,12pt,leqno]{amsart}

 \usepackage{amsmath}
\usepackage{amssymb}
\usepackage{amsthm}
\usepackage{graphicx}

\theoremstyle{plain} 
\newtheorem{theorem}{Theorem}
\newtheorem{proposition}{Proposition}
\newtheorem{remark}{Remark}

\begin{document}
\title[Measuring the mixing efficiency]{Measuring the mixing
  efficiency in a simple model of stirring:some analytical results and
  a quantitative study via Frequency Map Analysis} 

\author{Timoteo Carletti and Alessandro Margheri}

\address[T. Carletti]{Scuola Normale Superiore, piazza dei
  Cavalieri, 7, 56126 Pisa, Italy} 

\address[A. Margheri]{Fac. Ci\^encias de Lisboa and Centro de
  Matem\'atica e Aplica\c c\~oes Fundamentais, Av. Prof. Gama Pinto 2,
  1649-003 Lisboa, Portugal} 

\email[Timoteo Carletti]{t.carletti@sns.it}

\email[Alessandro Margheri]{margheri@ptmat.fc.ul.pt}

 

\begin{abstract}
We prove the existence of invariant 
curves for a $T$--periodic Hamiltonian system which models
 a fluid stirring in a cylindrical tank, when 
$T$ is small and the assigned stirring protocol is 
piecewise constant. Furthermore, using the Numerical Analysis of the
Fundamental Frequency of Laskar, we investigate numerically the break down 
of invariant curves as $T$ increases and we give a quantitative estimate
of the efficiency of the mixing.
\end{abstract}

\maketitle

\section{Introduction}

In  \cite{Ha}  it is studied  a simplified model of the stirring of an
ideal fluid in a cylindrical tank by an agitator. A Lagrangian
representation is considered for the motion,  which is assumed to be
completely two dimensional.  The tank thus degenerates to its
boundary circle, which has radius $R$. The agitator is modelled as a
point vortex of strength $\Gamma$ and its position inside the boundary
circle as a function of time, denoted by $z(t)$, is a prescribed
$T$--periodic function called {\em stirring protocol}. 
Introducing  the complex coordinate $\zeta=x+iy$,
 the motion of the fluid is governed by the following   nonautonomous
 $T$--periodic Hamiltonian system 
\begin{equation}\label{bate}
\dot{\overline\zeta}=\displaystyle\frac{\Gamma }{2\pi i}
\frac{|z(t)|^2-R^2}{(\zeta-z(t))(\zeta\overline{z(t)}-R^2)}\, ,
\end{equation}
whose corresponding  Hamilton function is given by
\begin{equation}
H(t,\zeta)=\displaystyle\frac{\Gamma}{2\pi}\ln\left | 
\frac{\zeta-z(t)}{\overline{z(t)}\zeta-R^2}\right | \, .
\end{equation}
System (\ref{bate}) is defined on the set $\{(t,\zeta)\in{\mathbb R}\times
{\mathbb C}: \zeta\neq z(t), \zeta\neq \frac{R^2}{\overline{z(t)}}\}$ and it
is considered for $\zeta$ belonging to the invariant disk
$D_R=\{\zeta\in{\mathbb C} : |\zeta|\leq R\}$.   
The following piecewise constant stirring protocol
\begin{equation}
\label{protocollo}
z(t)=\begin{cases}
+b & \quad t\in \Big [ nT, nT+T/2 \Big) \\
-b & \quad t\in \Big [ nT+T/2,(n+1)T \Big ) \quad n\in{\mathbb Z} \, ,
\end{cases}
\end{equation}
with fixed $0<b<R$, has been investigated in~\cite{Ha}. In this case 
  the motion can be integrated
over finite time, whereas the long time behavior of the  model is studied by
numerical experiments  when $b=1/2$  for different values of
$T.$ In ~\cite{Ha}, the objective of the author is  to study which
  features control the onset of chaos and hence the efficiency of the
  mixing.  
For this aim, the regimes of regular and chaotic behavior are identified and a  {\em qualitative}
  description (see~\cite{Ha} page 736 second paragraph) of the  mixing
  efficiency of the model has been provided  
 for several values of the parameters $(b,T)$. 
It turns out that the regions of chaotic behavior
consume a larger and larger portion of the phase space as $T$
increases, disrupting the regular pattern observed when $T$ is
very small (see Figure \ref{fig:pmap}). Quoting~\cite{Ha},
  '{\it...for small $T$ the
model [...] will look more and more like the two fixed agitator
system; [...] Thus one would expect 'convergence'  as $T\rightarrow 0$. }' 

In Section 2 of this work we will show that  averaging theory provides
the appropriate framework to investigate the above statement. Loosely
speaking, in such a setting we will be able to show that when $T\to 0$
the flow of system (\ref{bate}) with the stirring protocol
(\ref{protocollo}) converges in the $C^4$ topology  to the Hamiltonian
flow corresponding to two fixed agitators. As a consequence,  
by using  Moser's Small Twist Theorem (\cite{Mo},\cite{Or},\cite{He1},\cite{He2})
we will prove that  the regularity of the pattern observed in the
numerical experiments for small $T$ (see  Figure 2  a), b) and c) page
733 of~\cite{Ha} and Figure~\ref{fig:pmap} for $T=0.05$) is due, at
least in a  suitable  annular region inside $D_R$,  to the presence of
invariant curves. In fact, such curves  are, as one readily realizes,
an obstruction to the mixing of the fluid. 

The result which we present with this respect is the following:
\begin{center}
\begin{figure}[ht] 
 \begin{center}
  \makebox{\includegraphics[scale=0.15]{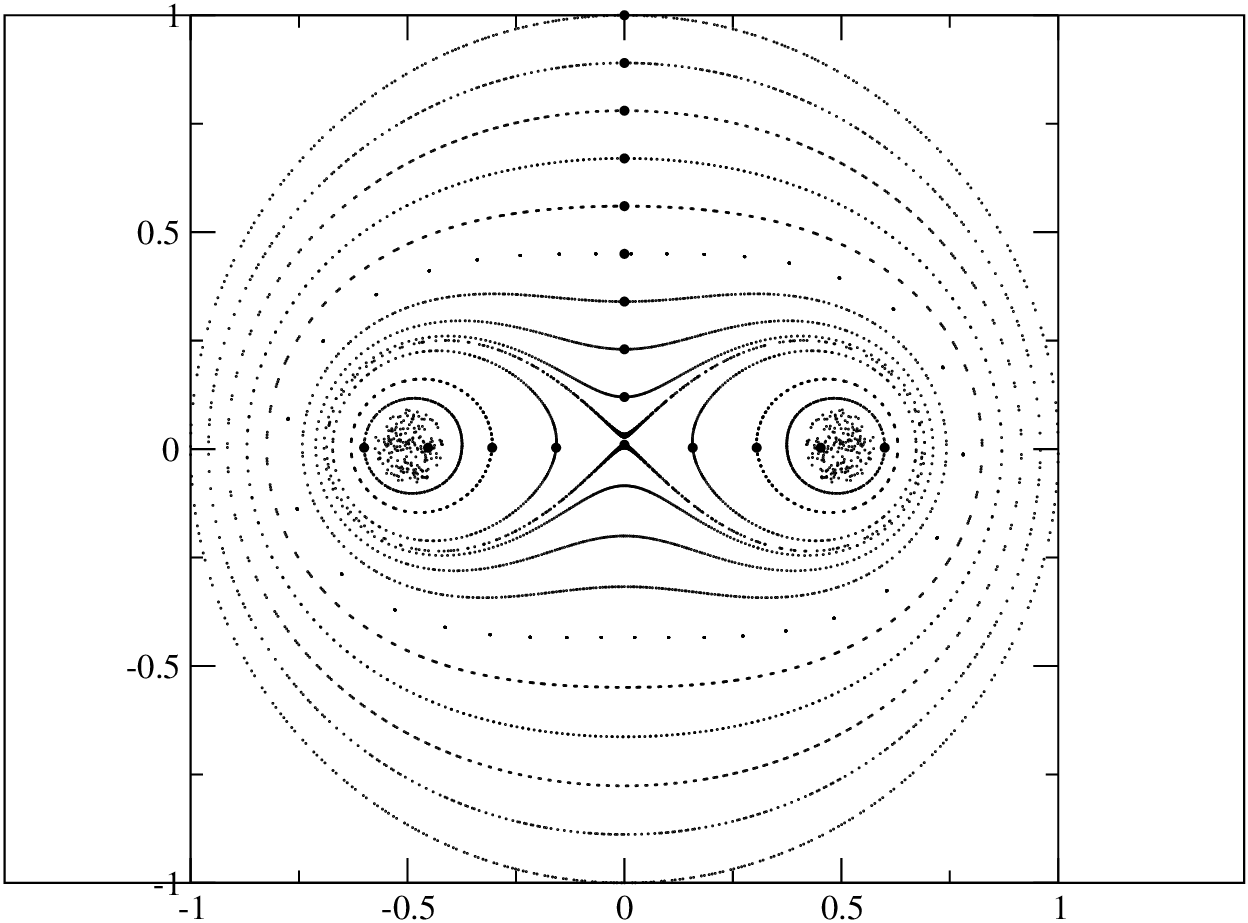}}\,
  \makebox{\includegraphics[scale=0.15]{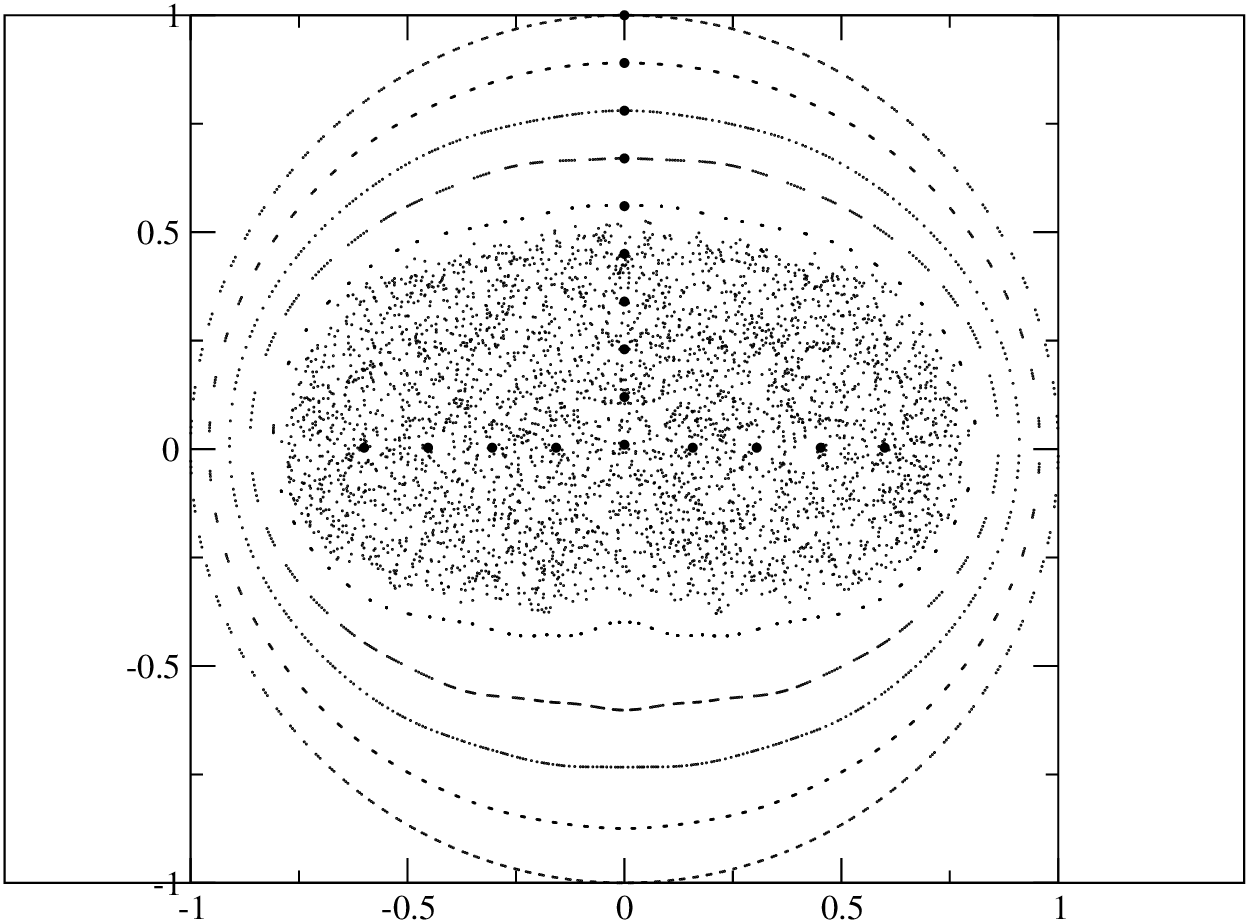}}\,
  \makebox{\includegraphics[scale=0.15]{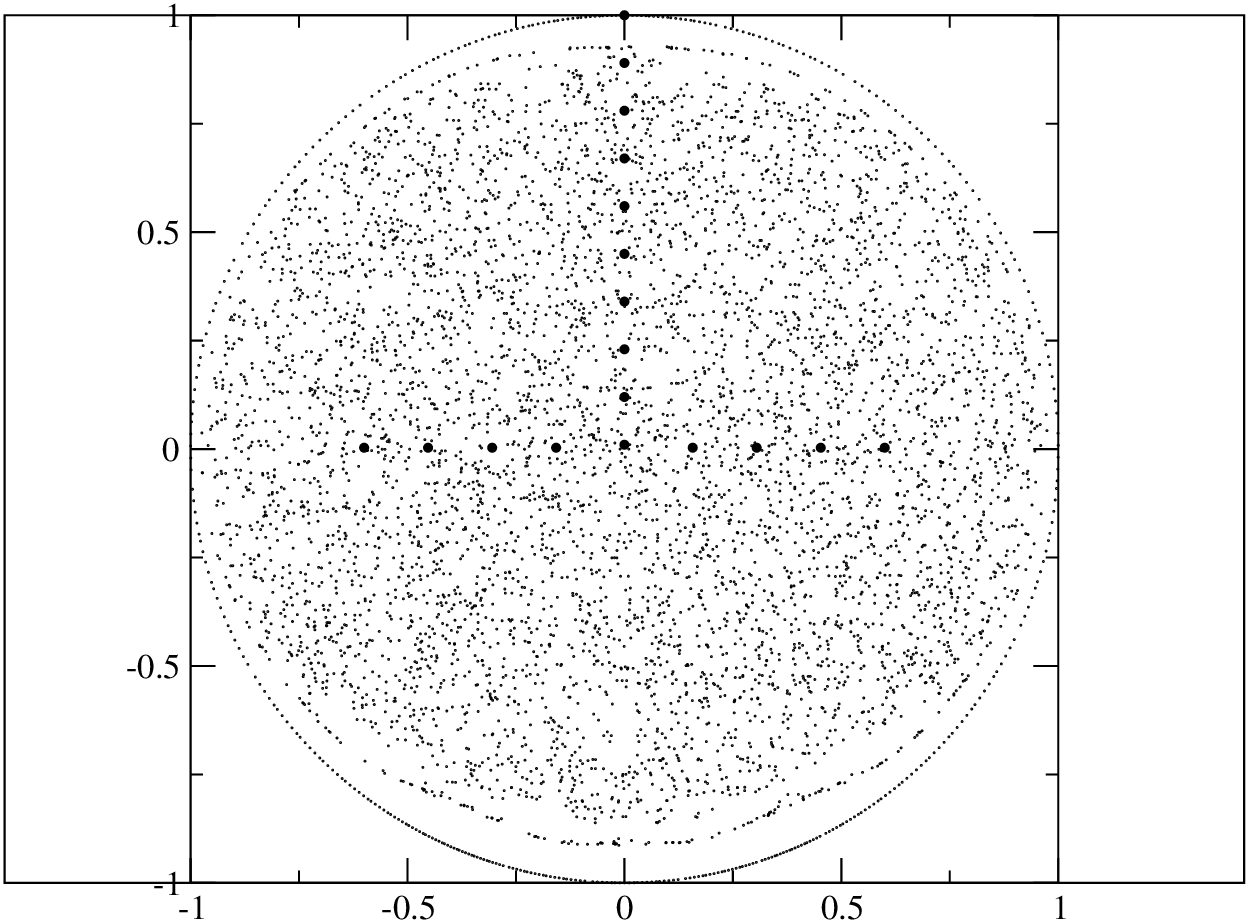}}
  \caption{Different behaviors of the Poincar\'e map, with $b=1/2$, 
$R=1$ and $T=0.05$ (left), $T=0.5$ (middle) and $T=1.5$ (right). Each plot is obtained following $10000$ time units the same $18$ initial data.}
  \label{fig:pmap}
  \end{center}
  \end{figure}
\end{center}

\begin{theorem}
\label{thm:main}
 Let $z(t)$ be  the stirring protocol (\ref{protocollo}) and let $T>0$ be its period. Then, there exists an annular region $A$ inside $D_R$ such that the $T$--Poincar\'e map associated to system~(\ref{bate})
 has invariant curves in $A$  for any  sufficiently small  $T$.
\end{theorem}

As $T$ increases the invariant curves break down. In Section 3, we will study numerically 
this phenomenon by means of the {\em Numerical Analysis of the
Fundamental Frequency, NAFF} of Laskar~\cite{Laskar}. This allows us to 
extend  our previous analytical result to cover larger $T$ ranges, by
showing numerical evidence that invariant curves persist close enough to
 the boundary of the invariant disk $D_R$ for large values
of $T$.
Furthermore, we are able to give a {\em quantitative} description of the
 efficiency of the mixing in function of the parameters $(b,T)$,
by measuring, 
through NAFF, the portion of phase space filled by invariant curves.

{\it Acknowledgements.}
Support from GRICES/CNR project is acknowledged. The second author
would also thank support from FCT

\section{Existence of invariant curves}
\label{sec:main}

In order to establish the existence of invariant curves for system
(\ref{bate}) with the piecewise constant  
protocol~(\ref{protocollo}) and small $T,$ we will show  that its
$T$--Poincar\'e operator (i.e. its period map) satisfies the
assumptions of  the so called {\em Moser's Small Twist Theorem}.  

Before outlining the strategy of our proof, it is worth to  recall the
statement of this important result. 

 Fixed  $I_1>I_0>0$,  let us consider the annulus defined by
 ${A}=\{(\bar{\theta}, I): \bar{\theta}\in {\mathbb S}^1\, , 
 I_0\leq I\leq I_1\}$.  Given a  mapping  $M:{A}\rightarrow {\mathbb R}^2$,
 we can find a lift  to the universal cover $\mathcal{A}= 
\{(\theta, I): \theta\in{\mathbb R}\, , I_0\leq I\leq I_1\}$ of $A$,  which we still denote by $M$. 

\medskip

\noindent{\bf Moser's Small Twist Theorem} {\it Let $\alpha$ be a
  $C^4([I_0,I_1])$ function satisfying 
\begin{equation}
\alpha'(I)<0 \quad \forall I\in [I_0,I_1] \, .
\end{equation}
Then, there exists $\epsilon>0$, depending on $I_1-I_0$ and $\alpha$, such 
that the map $M : A\mapsto {\mathbb R}^2$ has invariant curves if it satisfies 
the conditions below
\begin{itemize}
\item[a)] $M$ has the intersection property, that is,
 for any Jordan curve $\gamma$ homotopic to 
the circle $I=I_0$ in $A$,  $M(\gamma)\cap \gamma\not\equiv \emptyset$;
\item[b)] the lift of $M$ can be expressed in the form
\begin{equation}
M(\theta,I)=(\theta+T\alpha(I)+T\phi_1(\theta,I,T),I + T\phi_2(\theta,I,T))\, ,
\end{equation}
for some $T\in (0,1)$  and  $\phi_1, \phi_2\in C^4({\mathcal A})$ with  
$||\phi_1||_{C^4({\mathcal A})}+||\phi_2||_{C^4({\mathcal A})}<\epsilon$.
\end{itemize}}

This version of the theorem is presented in \cite{Or} and may be
proved using the techniques developed in \cite{He1,He2}.  

In what follows, a map of the form $(\theta+T\alpha(I),I )$ with
$\alpha'(I)\neq 0$, $I\in [I_0,I_1]$, will be referred to  as {\em small twist map}.

The application of  the previous theorem will be obtained performing
the following  steps. First,  by  using a classical construction from
averaging theory, in Subsection~\ref{ssec:averaged} we will rewrite  system
(\ref{bate})   as a $1$--periodic  perturbation of  the Hamiltonian
averaged system, being $T$ the small parameter. The averaged system
corresponds to the two fixed agitator model. Then, exploiting  the
geometry of the phase space of the averaged system, we will construct
explicitly the corresponding action variable on a suitable  subset  of
$D_R$. As a consequence,   we will be able to show  that on an
appropriate  annulus inside $D_R$ the  $1$--Poincar\'e map of the
averaged system is a small twist map. 

Next, the $C^4$ estimates needed to apply the Small Twist Theorem will
be  provided by a general result about the differentiability of a flow
with respect to parameters,  which we recall in Proposition \ref{Ort}
of Subsection~\ref{ssec:close}.

Finally, in Subsection~\ref{ssec:mainthm}, we collect all these fact
and complete the proof of Theorem 1.

\subsection{The averaged system}
\label{ssec:averaged}
In what follows, we  identify ${\mathbb C}$ with ${\mathbb R}^2$ and $\zeta=x+iy$ with $(x,y)$.   
 We rewrite the Hamilton equations of system~(\ref{bate}) in
 a compact real form as
\begin{equation}\label{origisis}
\dot\zeta = J\nabla_\zeta H(t,\zeta) \, ,
\end{equation}
where $J$ is the  standard $2\times 2$ symplectic matrix and
$\nabla_\zeta=(\frac{\partial}{\partial x},\frac{\partial}{\partial y})$.

The Hamilton function $H(t,\zeta)$ is piecewise autonomous and may be
 expressed in the form: 
\begin{equation*}
H(t,\zeta)= \phi(t)H_+(\zeta)+(1-\phi(t))H_-(\zeta)\, , 
\end{equation*} 
where
\begin{equation}
H_\pm(\zeta)=\displaystyle\frac{\Gamma }{2\pi}\ln\left |
 \frac{\zeta \mp b}{b\zeta\mp R^2}\right |\, ,
\end{equation}  
and $\phi(t)$ is the $T$--periodic extension of the restriction to 
$[0,T]$  of the characteristic function of $[0,T/2]$. Hence, defining 
$\mathcal{D}:={\mathbb R}^2\setminus \{\pm b,\pm R^2/b \}$,
 $H(t,\zeta)$ is
 smooth  on the set $({\mathbb R}\setminus\{k T/2, k\in{\mathbb Z} \})\times\mathcal{D}$.
 Moreover,
 $H$  and  its   derivatives  with respect to $\zeta$, which exist in
 $\mathcal{D}$ for any $t\in {\mathbb R}$, have  jump discontinuities at $t=k T/2$, $k\in {\mathbb Z}$. 

Next step is to consider $T$ as a small parameter and  to determine a
suitable comparison limit system when $T\rightarrow 0$. 
For this aim, as  in \cite{Ha}, we first rescale  time  in (\ref{origisis}) by $t=T\tau$,  
so normalizing  to $1$ the period of the stirring protocol and of the vector field. 
In the new time $\tau$,  setting   $H_1(\zeta,\tau):=H(\zeta,\tau T) $
and $y(\tau):=\zeta(\tau T)$,
 system (\ref{origisis}) takes the form 
\begin{equation}
\dot y(\tau) = T J\nabla_y  H_1(\tau,y)\, .
\label{bateT}
\end{equation}
To this system the  theory of averaging applies (see
\cite{GuHo}[Section 4.4]) as follows. 
 Denote by
\begin{equation}
\Hat{H_1}(y):=\int_0^1  H_1(\tau,y)\, d\tau\,=
\displaystyle\frac{\Gamma }{4\pi}
\ln\left | \frac{y^2-b^2}{y^2-R^4/b^2}\right |
\label{ham}
\end{equation}
and let  $V\subset{\mathbb C}$ be an open set such that its closure is
contained in $D$.  Then for small enough $T_0>0$  we can find a
symplectic~\footnote{Introducing canonical variables $y=(p,q)$ and
  $\eta=(P,Q)$ the requested transformation can be obtained trough the
 following generating function $S(q,P,\tau)=Pq+T\int_0^{\tau}\left[ \hat{H}_1(P,q)-H_1(P,q,s)\right]
  ds$. Thus $p=\partial_q S$ can be inverted if $T$ is small enough
  and trivially $S(q,P,\tau+1)=S(q,P,\tau)$ for all $(q,P,\tau)$ which
  implies~\eqref{condiniz}.},
close to identity, change  of coordinates of  the form   
\begin{equation}\label{nocan}
y=\eta+T w(\tau,\eta),\,\quad  (\tau,\eta,T)\in{\mathbb R}\times V\times
(0,T_0] \, ,
\end{equation} 
with $w$ 1--periodic in $\tau$ and such that 
\begin{equation}\label{condiniz}
w(0,\eta)=0\, ,
\end{equation}
which transforms system (\ref{bateT})   into system  
\begin{equation}\label{bateTw}
\dot \eta = TJ\nabla_\eta\hat{H_1}(\eta)+T^2 h(\tau,\eta, T) \, .
\end{equation}
By construction, the function $h(\tau,\eta,T)$  in  (\ref{bateTw}) is
$1$--periodic in $\tau$, has jump discontinuities at $\tau=k/2$,
$k\in{\mathbb Z}$ and it is smooth on $({\mathbb R}\setminus\{k/2, k\in{\mathbb Z} \})\times V\times[0,T_0]$.  

System (\ref{bateT})  is a perturbation of the  following (integrable)
{\em averaged system}
\begin{equation}
\label{systmed}
\dot \eta = TJ\nabla_{\eta}\Hat{H_1}(\eta) \,  .
\end{equation}

Henceforth, for simplicity, we will write $\hat H$ instead of $\hat H_1$.

We note that the Hamilton function $\Hat{H}$ corresponds to a two
point vortices system, one vortex  being located at $(+b,0)$ and the second at $(-b,0)$. 
Moreover, exploiting the geometry of the  phase space of system
(\ref{systmed}), we  can   construct explicitly the action variable
for this system outside the homoclinic loops surrounding the vortices
(see Figure~\ref {fig:pmap} on the left for $T=0.05$). This is done as follows.

By introducing symplectic polar coordinates $\zeta=\sqrt{2r}e^{i\psi}$
 and setting for notational convenience $E=e^{4\pi \Hat{H}/\Gamma}$ we can express 
the level lines of $\Hat{H}$ outside the homoclinic loops
\begin{equation}
r(\psi)=\frac{b^2R^4}{R^4+b^4}\cos 2\psi ,\quad \psi\in 
\left[-\frac{1}{4}\pi,\frac{1}{4}\pi\right]\cup
\left[\frac{3}{4}\pi,\frac{5}{4}\pi\right]\, ,
\end{equation}
in the form $r=r(\psi, E)$ with $\psi\in{\mathbb R}$ and $E\in (b^2/R^4, 1/R^2]$. 
We observe that $E=1/R^2$ corresponds  the boundary of
$D_R$, whereas $b^2/R^4$ corresponds to the homoclinic loops.

 Taking into account that $r(\psi,E)=r(-\psi,E)=r(\pi-\psi,E)$ (namely
the averaged system is invariant with respect to  $\zeta \rightarrow \bar\zeta$ and
 $\zeta \rightarrow -\bar\zeta$), the
 action variable $I$ is given by
\begin{equation}
\label{act}
I(E)=\frac{1}{2\pi}
\int_{
\Hat{H}(\sqrt{2r}e^{i\psi})=\frac{\Gamma}{4\pi}\ln E} r
\, d\psi=\frac{a(E)}{2\pi}\int_{0}^{\pi/2}\sqrt{b(E)
\cos^2 (2\psi)+c(E)}\,d\psi \, ,   
\end{equation}
where
\begin{gather*}
a(E):=\frac{1}{1-E^2b^4}\, ,\quad b(E):=4b^4(1-E^2R^4)^2 \\
c(E):=4(E^2b^4-1)(b^4-E^2R^8)\, . 
\end{gather*}

In next proposition, we use   the action variable constructed above to
show  that   the $1$--Poincar\'e operator  of the averaged system,
which henceforth we  
denote by $\hat M_1^T$, is a small twist map on an appropriate annular region inside $D_R$.
\begin{proposition}
\label{prop:main}
The set
\begin{equation*}
\mathcal{E}:= \Big \{ E\in \left(\frac{b^2}{R^4}, \frac{1}{R^2}\right] :
 I^{\prime}(E)\neq 0\text{ and } EI^{\prime\prime}(E)+I^{\prime}(E)\neq 0 \Big \}\, ,
\end{equation*}
is a non empty interval. Fix $E_0\in (\inf \mathcal{E}, \frac{1}{R^2})$ 
and consider the annulus 
\begin{equation}
\label{A}
A:=\Big\{\zeta\in D_R : \frac{\Gamma}{4\pi} \ln E_0\leq \Hat{H}(\zeta)\leq
 \frac{\Gamma}{4\pi} \ln \frac{1}{R^2}\Big \}\, .
\end{equation}
Then, $\hat M_1^T$ is a small twist map in $A$.
\end{proposition}

\proof
An easy computation shows that   $1/R^2\in \mathcal{E}.$ The rest of the statement about 
 $\mathcal{E}$ follows from the continuity of   the maps $E\mapsto
 I^{\prime}(E)$, $E\mapsto  EI^{\prime\prime}(E)+I^{\prime}(E)$.   
Let us consider  action--angle coordinates $(\theta,I)$  for the 
averaged Hamiltonian system  (\ref{systmed}), with $I$ given by
(\ref{act}). In the annular domain $A$, defined by~\eqref{A}, the map
$I=I(\hat H):=I(E(\hat H))$ has an inverse $\hat H=\hat H(I)$ on the
interval  
$[I_0,R^2/2]$, where $I_0:=I(E_0)$. 
 Moreover, on such interval we have
 $\hat H^{\prime\prime}(I)\neq 0$. Since, $\hat H^{\prime\prime}(1/R^2)<0$,
 it follows that 
\begin{equation}\label{twist}
\hat H^{\prime\prime}(I)<0,\quad  I\in[I_0,R^2/2] \, .
\end{equation} 
 The $1$--Poincar\'e map
 associated to system (\ref{systmed}) admits a lift to the set
 $\mathcal{A}:={\mathbb R}\times [I_0,R^2/2]$ given by  
 $\hat M_1^T(\theta,I)=(\theta+T \hat H^{\prime}(I),I)$ and therefore,
 taking also into account (\ref{twist}), it follows that it is a small
 twist map. 
\endproof

In order to apply Moser's Small Twist Theorem, we need some further 
estimates, which will be provided by the consequences of the
Peano's Theorem  stated below.

\subsection{Differentiability with respect to parameters}
\label{ssec:close}

Let $T>0$ be the period of the stirring protocol (\ref{protocollo})
and let $M_1^T$ be the $1$--Poincar\'e operator of system
(\ref{bateTw}). 
In  this subsection   we  will recall a general result, namely
Proposition \ref{Ort} below,  which is a consequence of the
differentiability of the flow of a system of O.D.E.'s with respect to
parameters. This proposition  is a restatement of  Proposition 6.4 in
\cite{Or}  under the  standard assumptions considered in dealing with
vector fields which are discontinuous in the independent variable. 
This result will provide the estimates needed to show  that, for
sufficiently small $T$,  the maps $M_1^T$ and $\hat{M}_1^T$ are
$C^4$--close  on  the annulus $A$ defined by (\ref{A}). 

This will be done  in the next subsection, and the proof of our main result will follow easily.

Let $V$ be an open subset of ${\mathbb R}^n,$  let $T_0$  be  a fixed positive
number and let $m$ be a nonnegative integer. 
Consider the following differential equation depending on one parameter  
\begin{equation}\label{equ}
\frac{dx}{dt} =F(t,x, T) \, ,
\end{equation}
where $ F:[0,1]\times V \times [0, T_0]\to {\mathbb R}^{n},\quad
(t,x,T)\mapsto F(t,x,T)$  is a function which  satisfies 
\begin{itemize}
\item the map $(x,T)\mapsto F(t,x,T)\in 
C^{m+1}(V\times [0,T_0])\,$  for almost all $t\in [0,1]$;
\item for $0\leq j+k\leq {m+1}$, the maps $t\mapsto \frac{\partial^{k+j}}
{\partial x^k\partial T^j }F(t,x,T)$ are $t$ 
measurable for any $(x,T)\in V\times [0,T_0]$;
\item  for each compact set $K\subset V$  
there exists $C>0$ such that $$\Big |\frac{\partial^{k+j}}
{\partial x^k\partial T^j} F(t,x,T)\Big |\leq C \, ,$$ 
for all
$0\leq k+j \leq m+1$ and $(t,x,T)\in [0,1]\times K\times
[0,T_0]$.
\end{itemize}

The solution of (\ref{equ}) satisfying $x(0)= x_0$, will be denoted by 
$x(t; x_0,T)$.  By the general theory of ordinary differential
equations, $x$ is of class $C^{0,m+1,m+1}$ in its three arguments
whenever it is defined. The following result is a consequence of this
fact. 

\begin{proposition}\label{Ort}
 Let $A$   be a compact subset of $V$ such that  for every $x_0\in A$
 and $T\in [0,T_0]$ the solution $x(t;x_0,T)$ is well defined in
 $[0,1]$. Then, for each $(t,x_0,T)\in [0,1]\times A\times [0,T_0]$
 the expansion below holds 
$$x(t;x_0,T)=x(t;x_0,0)+\frac{\partial x}{\partial
  T}(t;x_0,0)T+R(t;x_0,T)T\, ,$$
where the remainder  $R$ satisfies
$$||R(t;\cdot,T)||_{C^m(A)}\to 0,\quad T\to 0 \, ,$$
uniformly in $ t\in [0,1]$.
\end{proposition}

\subsection{Proof of Theorem \ref{thm:main}}
\label{ssec:mainthm}

We are now able to prove our main theorem. We will show  that for sufficiently small $T$ the
 $T$--Poincar\'e map associated to system (\ref{bate}) with the
 stirring protocol (\ref{protocollo}) has invariant  
 curves in the set $A$ defined by (\ref{A}).

For a fixed $R^{\prime}$ satisfying  $R<R^{\prime}<R^2/b$, we define
\begin{equation} 
V:= \Big \{ \zeta=\sqrt{2r}e^{i\psi}\in{\mathbb C} : \frac{b^2R^4}{R^4+b^4}\cos 2\psi
 <r<(R^{\prime})^2/2 \Big \}\, .
\end{equation}

This set, which contains  $A$, corresponds   to the open disk of radius $R^{\prime}$ minus the compact 
 region bounded by the homoclinic loops. Hence, it has
 positive distance from the singularities of $\Hat{H}$.
 Now we choose a sufficiently small $T_0>0$  for which  the
 transformation (\ref{nocan}) is well defined on $[0,1]\times V\times
 [0,T_0]$  and, moreover, all the solutions of system (\ref{bateT})
 with initial conditions in $A$ are well defined  in    
$[0,1]$. 
 Let $T<T_0$  be the period of $z(t).$ 
 By Proposition \ref{Ort}, and by introducing the action angle
 coordinates considered in Proposition  \ref{prop:main}, we can find a
 lift of  $M_1^T$ to $\mathcal{A}:={\mathbb R}\times [I_0,R^2/2]$ of the form 
\begin{equation*}
M_1^T(\theta,I)=(\theta + \hat{H}^{\prime}(I)T,I)+R(\theta,I,T)T \, ,
\end{equation*}
  where $||R(\cdot,T)||_{C^4(\mathcal{A})}\rightarrow 0$  for $T\rightarrow 0$.
Hence,  $M^T_1$ satisfies hypothesis b) of Moser's
Small Twist Theorem for $T$ small enough.
 Moreover, $M_1^T$ has the intersection property in $A$, being an
 area--preserving map in $A$ for which the boundary of $A$ is
 invariant.

Then, by Moser's Small Twist Theorem, $M_1^T$ has many invariant curves in $A$ for  $T$
 small enough, and so does the $T$--Poincar\'e operator of system
 (\ref{bate}), which  obviously  coincide with 
$M_1^T$.
\endproof

\section{Numerical Results}
\label{sec:numres}

In this section we will study system   (\ref{bate}) with the stirring
protocol (\ref{protocollo}) from a numerical point of view. More
precisely, we will be interested in two types of numerical
experiments. In the first one, we will explore the 
{\em continuation} properties with respect to $T$  of the invariant curves of the
$T$--Poincar\'e map of system (\ref{bate}) obtained for small $T$ in Theorem \ref{thm:main}.
 In \S \ref{applicI},  we will give  numerical evidence of their
 persistence close to the boundary for large values of $T$. 

In the second experiment,  we will be able to give a {\em quantitative
 }  description  of the  mixing efficiency  
 of the piecewise constant stirring protocol. This will be done in  \S
 \ref{applicII}, where we will  investigate the parameters $(b,T)$
 giving rise to an efficient stirring.  

Both these numerical experiments will be carried out by 
using  the {\em Numerical Analysis of the Fundamental Frequency (NAFF)}
of Laskar\cite{Laskar}.

 This is a numerical method which allows to
obtain a global view of the behavior of a dynamical system 
by studying the properties of the frequency map, numerically defined 
from the action--like variables to the frequency space using adapted 
Fourier techniques. This method has been used to investigate  a wide 
class of dynamical systems, such as  the solar system~\cite{Laskar2},
the galactic dynamics~\cite{PAP96}, particle 
accelerators~\cite{LaskarRobin} and  the
 standard map~\cite{CarlettiLaskar}.

For  sake of completeness, let us first present briefly the main outlines
 of the method, referring to~\cite{Laskar} for further details.

 \subsection{Frequency Map Analysis.}

 Let us consider an  $n$--degrees of freedom quasi--integrable 
 Hamiltonian system, expressed in  action--angle variables by 
\begin{equation}
  H\left(I,\theta;\epsilon\right)=H_0\left(I\right)+
\epsilon H_1\left(I,\theta\right) \, ,
\label{eq:hamsys}
\end{equation}
where $H$ is a real analytic function  of $\left(I,\theta\right)\in B\times {\mathbb T}^n$, 
 $B$ is an open domain in 
${\mathbb R}^n$, $\,{\mathbb T}^n$ is the $n$--dimensional torus and $\epsilon$ is a  small real parameter.
For $\epsilon=0$ the system is integrable:
the motion takes place on invariant tori $I_j=I_j(0)$ described at constant
velocity $\nu_j\left(I\right)=\frac{\partial H_0}{\partial I_j}
\Big\rvert_{I\left(0\right)}$, for $j=1,\ldots,n$.
Assuming a non--degeneration condition on $H_0,$ the frequency map 
$F:B\rightarrow {\mathbb R}^n$
\begin{equation*}
F:I\mapsto F\left(I\right)=\nu \, , 
\end{equation*}
is a diffeomorphism onto its image $\Omega$. In this case KAM 
theory~\cite{Kolmogorov,Arnold,Mo}
ensures that for sufficiently small values of $|\epsilon|$, there exists a 
Cantor set
$\Omega_{\epsilon}\subset \Omega$ of frequency vectors satisfying 
a Diophantine condition, for which the 
quasi--integrable
system with Hamilton function~(\ref{eq:hamsys}) still possess smooth invariant tori. These tori are 
$\epsilon$--close to those of the
unperturbed system, and support the linear flow $t\mapsto \theta_j\left(t\right)=
\nu_j t+\theta_j\left(0\right) (\mod 2\pi)$ for
$j=1,\ldots,n$. Moreover, according to P\"oschel~\cite{poschel} there 
exists a diffeomorphism
$\Psi:{\mathbb T}^n\times \Omega \rightarrow {\mathbb T}^n\times B$
\begin{equation*}
  \Psi:\left(\phi,\nu\right)\mapsto\left(\theta,I\right)\, ,
\end{equation*}
which is analytic with respect to $\phi$ in ${\mathbb T}^n$ and 
$\mathcal{C}^{\infty}$ w.r.t $\nu$ in
$\Omega_{\epsilon}$, and which transforms the Hamiltonian 
system generated by~(\ref{eq:hamsys}) into
\begin{equation*}
  \begin{cases}
   \frac{d\nu_j}{dt}\left(t\right)=0\\ 
   \frac{d\phi_j}{dt}\left(t\right)=\nu_j \, .
  \end{cases}
\end{equation*}
For frequency vectors $\nu \in \Omega_{\epsilon},$ the invariant torus 
can be represented in the
complex variables $(z_j=I_je^{i\theta_j})_{j=1,n}$  by a quasi--periodic 
function
\begin{equation}
  \label{eq:cmplxtorus}
  z_j\left(t\right)=z_j\left(0\right)e^{i\nu_jt}+\sum_{m}a_{j,m}
\left(\nu\right)e^{i <m,\nu> t} \, .
\end{equation}
Taking a section $\theta=\theta_0$ of the phase space, for some 
$\theta_0 \in {\mathbb T}^n$, we
obtain the frequency map $F_{\theta_0}:B\rightarrow \Omega$ 
\begin{equation}
  \label{eq:freqmap}
  F_{\theta_0}:I\mapsto \pi_2\left(\Psi^{-1}\left(\theta_0,I\right)\right) \, ,
\end{equation}
where $\pi_2\left(\phi,\nu\right)=\nu$ is the projection onto $\Omega$. 
For sufficiently small
$|\epsilon|$ the non--degeneration condition ensures that $F_{\theta_0}$ is a 
smooth diffeomorphism.

If we have a numerical (complex) signal over a finite time
 span $\left[ -K, K\right]$ and we want to recover a quasi--periodic 
structure, we can construct an 
$N$--terms quasi--periodic approximation 
$\hat f \left( t \right)=\sum_{k=1}^N 
a^{(K)}_k e^{it\nu^{(K)}_k }$. Frequencies 
$\nu^{(K)}_k$ and amplitudes $a^{(K)}_k$ are determined with
an iterative scheme (possibly) involving some weight function 
(Hanning filter). 

Assuming some ''good'' arithmetic properties of the frequencies and 
using the $p$--th Hanning filter, then one 
can prove \cite{Laskar} that NAFF
 converges towards the first\footnote{We assume to enumerate frequencies 
according to decreasing amplitudes: $|a_k^{(K)}|\geq |a_{k+1}^{(K)}|$. That
is, some smoothness of the signal is assumed.} ''true''
frequency of the given signal, $\nu_1$, with the following  asymptotic
expression for $K \rightarrow +\infty$  
\begin{equation}
\label{eq:asintexpr}
|\nu_1 - \nu_1^{(K)} |= \mathcal{O}\! \left(\frac{1}{K^{2p+2}}\right) \, ,
\end{equation}
which will usually be several
order of magnitude better than the $\mathcal{O}\left(1/K\right)$ order obtained with FFT.

Using the NAFF algorithm, it is  possible to construct numerically a 
frequency map (see 
Figure \ref{fig:freqmap11}) in the following way:
\begin{itemize}
\item[a)] fix all angles to some value $\theta_{0}$;
\item[b)] for all initial values $I_0$ of the action variables,  integrate 
numerically the trajectories with initial condition 
$(I_0,\theta_0)$ over the time span $K$;
\item[c)]  look for a quasi--periodic approximation of the trajectory,  
with the previous algorithm identifying the fundamental 
frequency $\nu_1^{(K)}$ of this 
quasi--periodic approximation. 
\end{itemize}

Because of (\ref{eq:asintexpr}), we can use NAFF algorithm to test the 
goodness of the reconstructed quasi--periodic
 approximation of the given signal, showing numerical evidence that motion
 doesn't take place on invariant tori because of some ''diffusion in frequency
space''. More precisely, this can be done by replacing step c) of the
previous algorithm with the following  

\begin{itemize}
\label{pag:11}
\item[c')] divide the time span $[-K,K]$ into smaller parts (possibly
overlapping): $[t_l,t_l+K_1]$, for some $t_l$. Then NAFF reconstructs for 
each $l$ a frequency $\nu^{(K_1)}(l)$. If the frequencies 
$(\nu^{(K_1)}(l))_l$ 
coincide, up to some prescribed numerical precision, when $l$ varies, then
we can conclude that we obtained a good quasi--periodic
approximation. Otherwise, we can measure the diffusion in the
frequency space. 
\end{itemize}

\subsection{Application I: existence of invariant curves}
\label{applicI}

Let us  consider the problem of the existence of invariant curves
homotopic to the boundary of $D_R$  for the 
$T$--Poincar\'e map of system (\ref{bate})  with stirring protocol (\ref{protocollo}). One 
readily realizes that every such curve 
will intersect the vertical segment $\ell:=\{ \zeta\in {\mathbb C} : x=0 \, ,
\, y \in [0,R]\}$.

In order to  reconstruct a frequency map, we fix a piecewise constant stirring protocol
by choosing the parameters $(b,T).$ Then, for $N$ initial data on the segment 
$\ell$ we construct numerically the $T$--Poincar\'e map and we iterate it 
$N_{iter}$ times.

\begin{remark}[Numerical computation of the Poincar\'e map]
\label{rem:integration}
 As already remarked in \cite{Ha}, one can construct the
 Poincar\'e map for the piecewise constant protocol avoiding the 
integration of the Hamiltonian system (\ref{bate}), which results
CPU--time consuming and introduces additional errors. The idea is that in
each half period the Hamiltonian system is autonomous, hence integrable and,
moreover, the integration can be explicitly done: for $t\in[0,T/2)$, 
the motion takes place
on an arc of circle with radius $\rho = \frac{\lambda}{1-\lambda^2}\left(
\frac{R^2}{b}-b\right)$, centered at $\zeta_{c}=\frac{b^2-\lambda^2 R^2}
{b\left(1-\lambda^2\right)}$, with $\lambda = b\Big \lvert
\frac{\zeta(0)-b}{b\zeta(0)-R^2}\Big \rvert$. Thus, setting $\zeta(0)=\zeta_c +
\rho e^{i\phi(0)}$, after half period the point will be at $\zeta(T/2)=
\zeta_c +\rho e^{i\phi(T/2)}$ where: 
\begin{equation}
\label{eq:mezzamappa}
\phi(T/2)-\frac{2\lambda}
{1+\lambda^2}\sin \phi(T/2) = \frac{\Gamma}{2\pi \rho^2}\frac{1-\lambda^2}
{1+\lambda^2}\frac{T}{2}\, .
\end{equation}
The motion in the second half  period is similar.

Equation (\ref{eq:mezzamappa}) can be solved using a sixth order
Newton's--type method. However, when
$b$ is close to $R$ and/or $T$ is large, Newton's algorithm doesn't
converge to the good solution unless we provide a very accurate  approximation of $\phi(T/2)$ as 
initial value for the algorithm. In order  to overcome this difficulty, we used the 
following trick. Equation~(\ref{eq:mezzamappa}) is nothing but {\em
  Kepler's equation} and in \cite{watson} an explicit solution is given
in terms of Bessel's functions (which results, in computations, less efficient than
our Newton's--type method). Neverthless for large values of the parameters we obtain a good 
approximation to start with  Newton's algorithm by summing few terms of the
previous explicit solution.
\end{remark}

\begin{center}
  \begin{figure}[h] 
  \begin{center}
   \makebox{\includegraphics[scale=0.35]{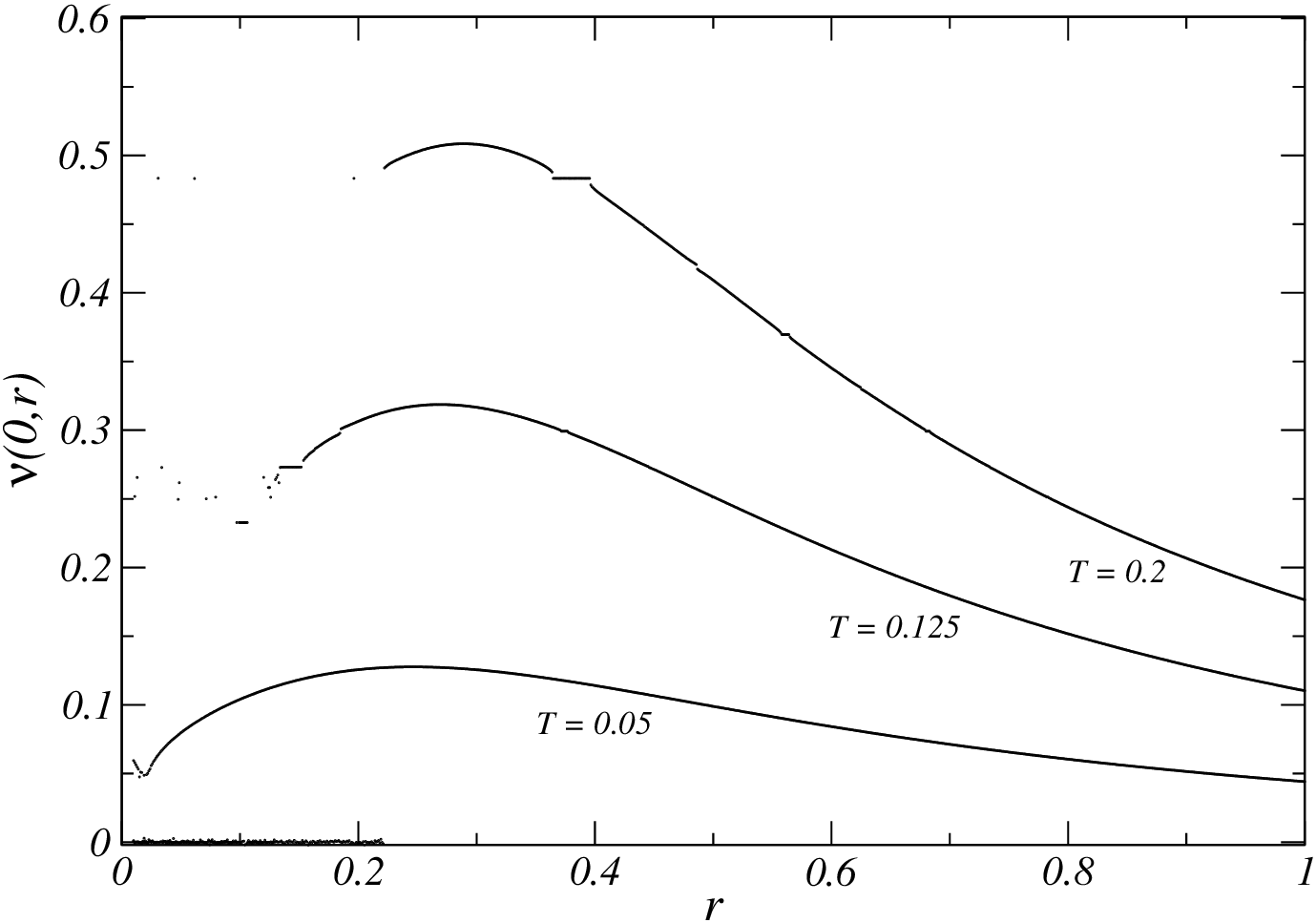}}
  \caption{The Frequency Map numerically reconstructed by NAFF. Parameters
 are $b=1/2$, $R=1$ and $T\in \{0.05,0.125,0.2 \}$. Initial data
 are $1000$ points $\zeta_j(0)=ir_j$, equally spaced on $[0,i]$ and orbits are
computed for $N_{iter}=50000$.}
  \label{fig:freqmap11}
  \end{center}
  \end{figure}
\end{center}

\begin{center}
  \begin{figure}[h] 
  \begin{center}
   \makebox{\includegraphics[scale=0.35]{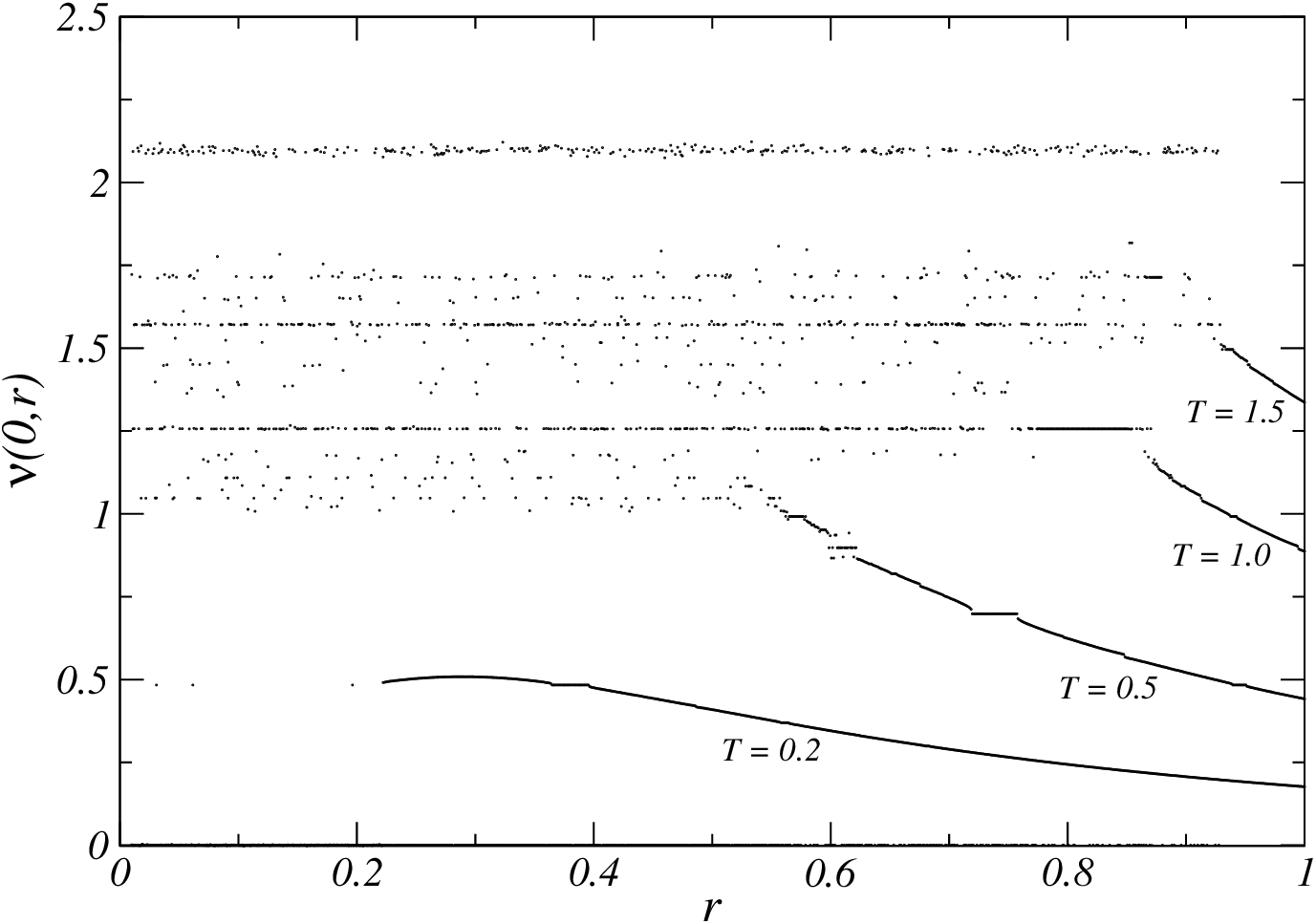}}
  \caption{The Frequency Map numerically reconstructed by NAFF. Parameters
 are $b=1/2$, $R=1$ and $T\in \{0.2,0.5,1.0,1.5 \}$. Initial data
 are $1000$ points $\zeta_j(0)=ir_j$, equally spaced on $[0,i]$ and orbits are
computed for $N_{iter}=50000$.}
  \label{fig:freqmap12}
  \end{center}
  \end{figure}
\end{center}

Once we have an orbit, $\left(M^k_T(\zeta(0))\right)_{0\leq k \leq N_{iter}}$,
 we use the NAFF algorithm to reconstruct a 
quasi--periodic approximation of it. The frequency map is not smooth at all  
because of the presence of resonances: there is no invariant curve with
rational rotation number. Thus, for a fixed resolution, the existence of
invariant curves is assumed if the frequency map looks ''regular'' (see
 Figure \ref{fig:freqmap11} for $T=0.05$), whereas
to ''irregular'' graphs we associate non existence of invariant curves
(see Figure \ref{fig:freqmap12} for $T\in \{0.5,1.0,1.5\}$). 

In Figure \ref{fig:freqmap11} we present some numerical results for the 
piecewise constant protocol (\ref{protocollo}) with $b=1/2$ 
and several values of $T$. On the left, for very small $T$, say $0.05$,
the curve looks regular, hence a large part of the disk is occupied by invariant
curves. We can observe that the map {\em is not twist} in the whole disk: there
is a point where the derivative of the frequency map is zero. For slightly 
larger $T$, say  $T\in(0.125,0.2)$, the curve is still
regular,  but an irregular pattern is showed close to the origin. Roughly 
speaking, no invariant curves are present for $T=0.2$ in the disk 
$|\zeta | \leq 0.22$. Also, elliptic points due to resonances are showed.

In Figure \ref{fig:freqmap12}, when $T$ increases, the irregular patterns consume larger and 
larger portions of the invariant disk, but still a regular frequency map is
present close to the boundary of $D_R$.

\begin{center}
  \begin{figure}[ht] 
  \begin{center}
   \makebox{\includegraphics[scale=0.35]{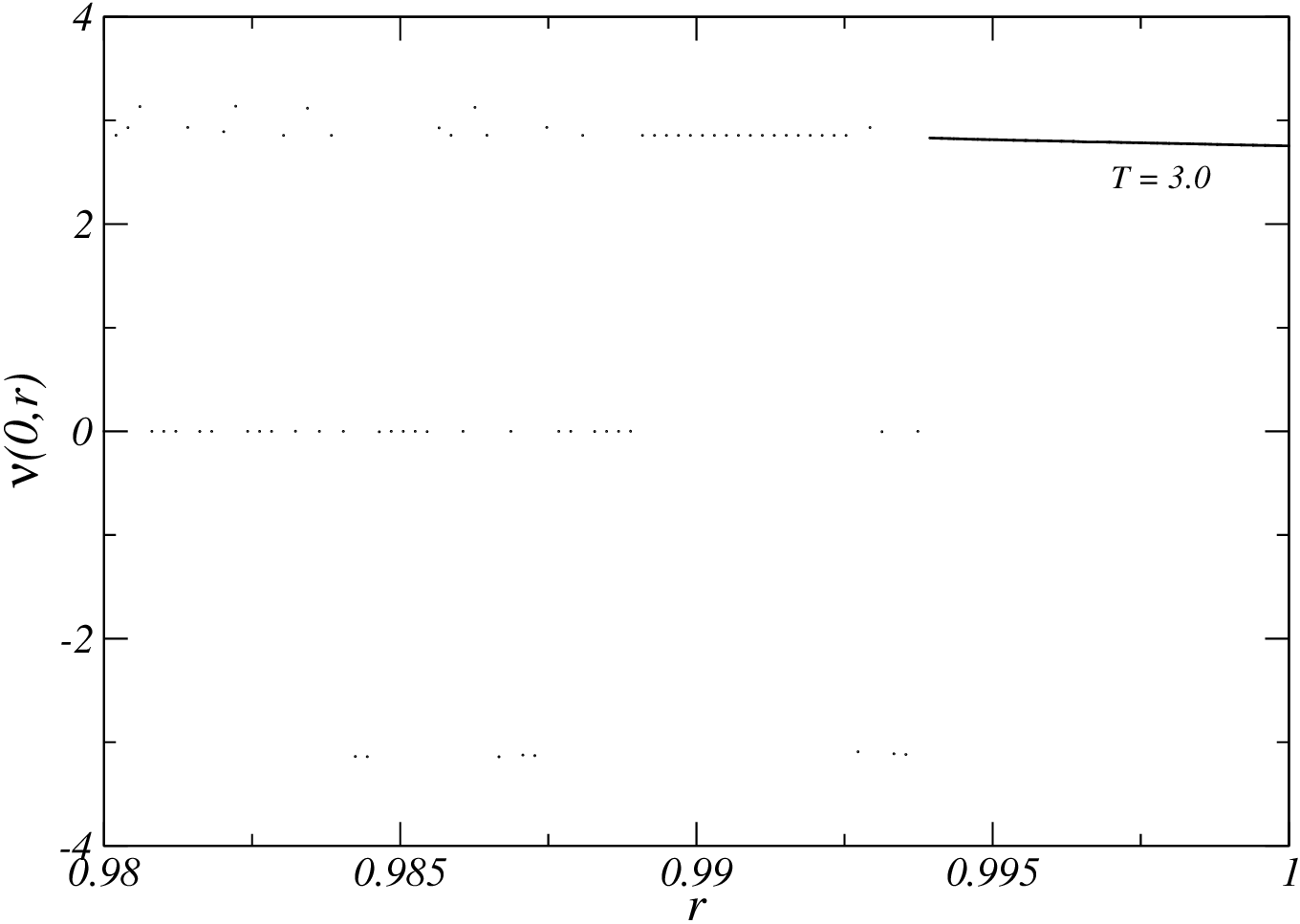}}
  \caption{The Frequency Map numerically reconstructed by NAFF. Parameters
 are $b=1/2$, $R=1$ and $T=3.0$.  Initial data
 are $1000$ points $\zeta_j(0)=ir_j$, equally spaced on $[0.98i,i]$ and
 orbits are 
computed for $N_{iter}=500000$.}
  \label{fig:freqmap2}
  \end{center}
  \end{figure}
\end{center}

This behavior of the system is displayed in Figure \ref{fig:freqmap2}, where we show a numerical result for $b=1/2$ and 
$T=3.0$, a period which is $60$ times larger than the smallest value in Figure 
\ref{fig:freqmap11}.  Note the change in the scale on
the axes w.r.t. to Figure \ref{fig:freqmap11} and \ref{fig:freqmap12} (the scale on the vertical
axis has been  changed to  show the large oscillations of the frequency map).

Other results for larger values of $T$ exhibit a similar behavior. This supports the conjecture that 
 invariant curves persist, close to the boundary, for arbitrarily large values of  $T$.

\subsection{Application II: mixing efficiency}
\label{applicII}

We have already observed that 
invariant curves are an obstruction to global mixing. 
We note now  that close to {\em robust} invariant curves (i.e. associate to good
frequencies) there is a neighborhood filled by invariant curves 
\cite{morbidelligiorgilli,morbidelligiorgilli1}, which prevents also
from  local mixing. Using  
the precision of NAFF and observation c') in page \pageref{pag:11}, we 
can test the efficiency of the stirring protocol by evaluating the goodness
of the reconstructed signal.

More precisely, we choose a couple of parameters $(b,T)$, then we divide 
the invariant disk $D_R$ into a fine grid of 
$N_{tot}$ points. For each point we determine numerically 
(see Remark \ref{rem:integration}) the orbit for a time interval $[0,K]$. 
Then, fixing some $K_1 < K$ and some positive $(t_l)_l$ we divide the orbit
into pieces $\left(\zeta(t)\right)_{t\in [t_l,t_l+K_1]}$ such that
intervals $[t_l,t_l+K_1]$ overlap. On each piece, the  
use of NAFF  gives us a fundamental frequency $\nu^{(K_1)}(l)$, where we
emphasize the dependence of the frequency on the $l$--piece of orbit.

For a given orbit with initial datum $\zeta_0=(x_0,y_0)$, let us define
\begin{equation}
\label{eq:indicatore}
\varepsilon(x_0,y_0) = -\log \Big\lvert 2 \frac{\max_l \nu^{(K_1)}(l) -
 \min_l\nu^{(K_1)}(l)}{\max_l \nu^{(K_1)}(l) +
 \min_l\nu^{(K_1)}(l)}\Big\rvert \, .
\end{equation}
This is a good indicator of the robustness of the orbit. In fact, if
$\varepsilon(x_0,y_0)$ is ''large'' then $\max_l \nu^{(K_1)}(l)$ and 
$\min_l\nu^{(K_1)}(l)$ are ''close together''. Thus, the reconstructed 
frequencies on different pieces of orbit  do not vary too much, and we can
assume that we are on a quasi--periodic orbit. On the other hand, if 
$\varepsilon(x_0,y_0)$ is ''small'' then $\max_l \nu^{(K_1)}(l)$ and 
$\min_l\nu^{(K_1)}(l)$ are ''far from each other''. In this case, since the reconstructed frequency is not
constant,  we conclude that we are not on a quasi--periodic orbit.
 Intermediate values of $\varepsilon(x_0,y_0)$ give rise to intermediate degrees of robustness.

\begin{center}
  \begin{figure}[ht] 
   \makebox{\includegraphics[scale=0.23]{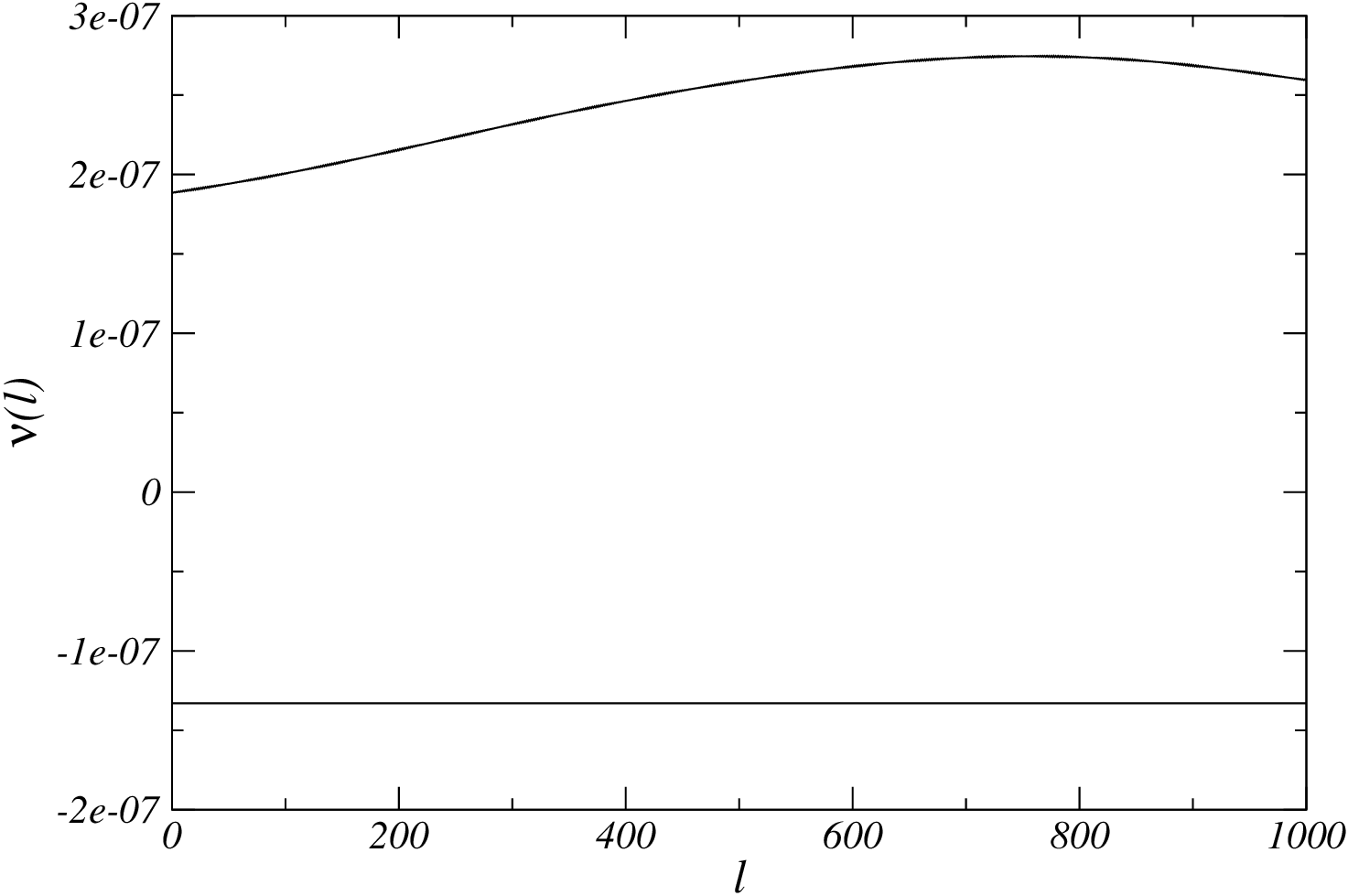}}\,\,\,\,
   \makebox{\includegraphics[scale=0.23]{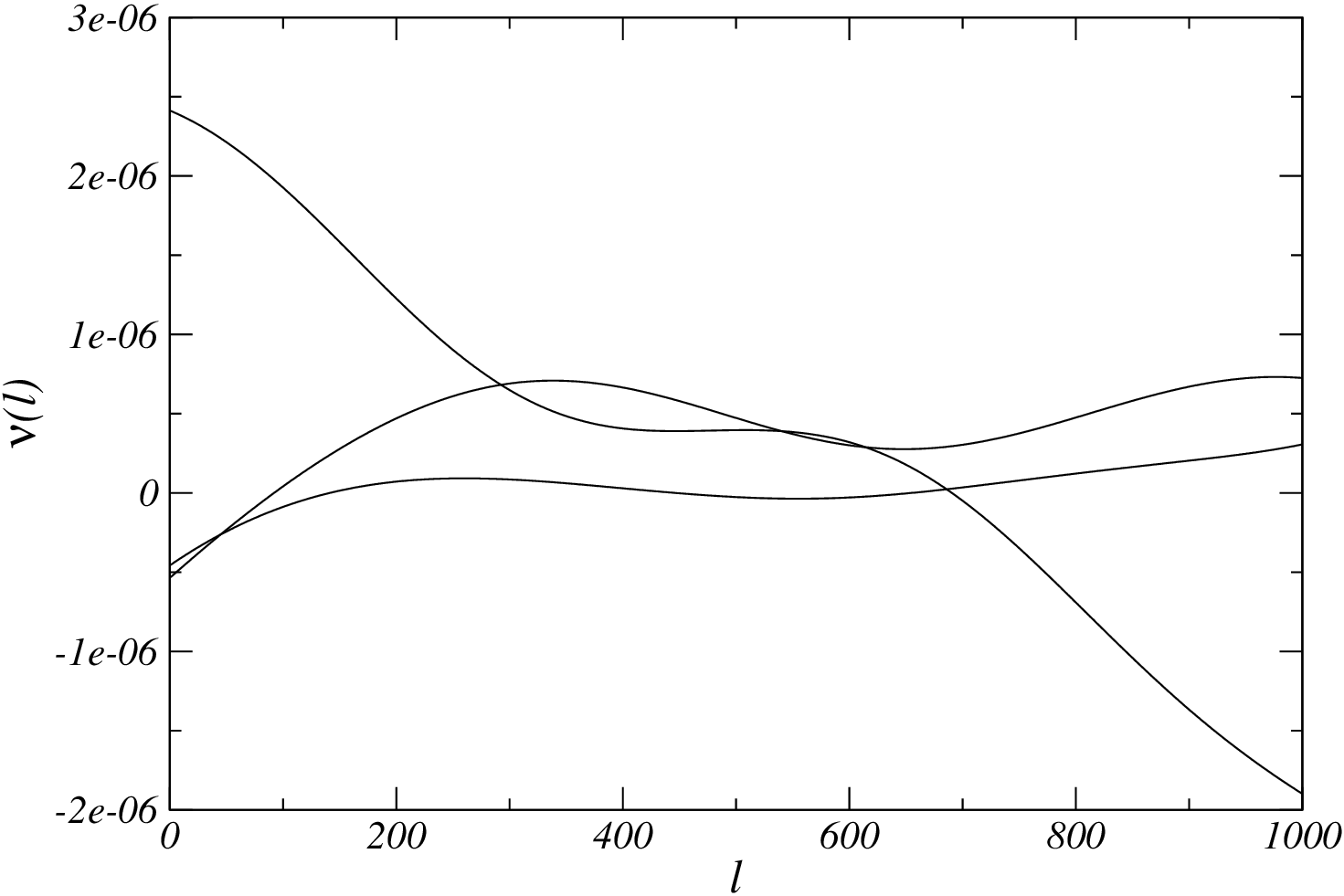}}
  \caption{Distance between the frequency $\nu_1^{(K_1)}(l)$,
   determined by NAFF, and the mean value of the $\nu_1^{(K_1)}(l)$'s,
   in function of the part of orbit $[t_l,t_l+K_1]$. On the left very
   regular orbits, on the right less regular ones. Both plots are made    with $N_{iter}=100000$, $K_1 = 50000$ and $t_l = 50(l-1)$.}
  \label{fig:rob}
  \end{figure}
\end{center}

We now fix a threshold value $\varepsilon_{thr}>0$ and we measure the portion
of initial data in $D_R$ to which one can associate a 
$\varepsilon_{thr}$--robust orbit
\begin{equation}
\label{eq:mdef}
m_{\varepsilon_{thr}}=\frac{\#\{ (x_0,y_0): \varepsilon(x_0,y_0) >
 \varepsilon_{thr} \}}{\#\{ (x_0,y_0)\}}\, .
\end{equation}
A preliminary analysis, see Figure~\ref{fig:rob}, of some robust orbits
allows us to choose an appropriate threshold value.
In this way we are able to turn the {\em qualitative} analysis of Figure 7
 page 738 of \cite{Ha} into a {\em quantitative} one.

With this respect, we present two results, the first with $\varepsilon_{thr}=12$, corresponding 
to $\Big\lvert \frac{\max_l \nu^{(T_1)}(l)}
{\min_l\nu^{(T_1)}(l)}-1\Big\rvert \leq 3\, 10^{-6}$,
 and the second with $\varepsilon_{thr}=15$ and $\Big\lvert 
\frac{\max_l \nu^{(T_1)}(l)}{\min_l\nu^{(T_1)}(l)}-1\Big\rvert \leq 1.5
\, 10^{-7}$. Each result has been obtained dividing the invariant disk
in approximately $30000$ points, equally spaced in both $x$ and $y$ by
$0.01$. Then, several values of $b$ and $T$ have been considered.

Finally, by using the previous ideas, we  precise as follows   the classification scheme given for a regime  
 in \cite{Ha}:
\begin{itemize}
\item[{\bf [I]}] {\em integrable} \,\, if $0.6 < m_{\varepsilon_{thr}} \leq 1$ (very poor mixing property)
;
\item[{\bf [T]}] {\em transitional} \,\, if $0.3 < m_{\varepsilon_{thr}} \leq 0.6$;
\item[{\bf [C]}] {\em chaotic} \,\, if $0< m_{\varepsilon_{thr}} \leq 0.3$  (efficient mixing),
\end{itemize}
and we  summarize our results in  the following Figure \ref{fig:effic}
\begin{center}
 \begin{figure}[ht] 
  \begin{center}
   \makebox{\includegraphics[scale=0.23]{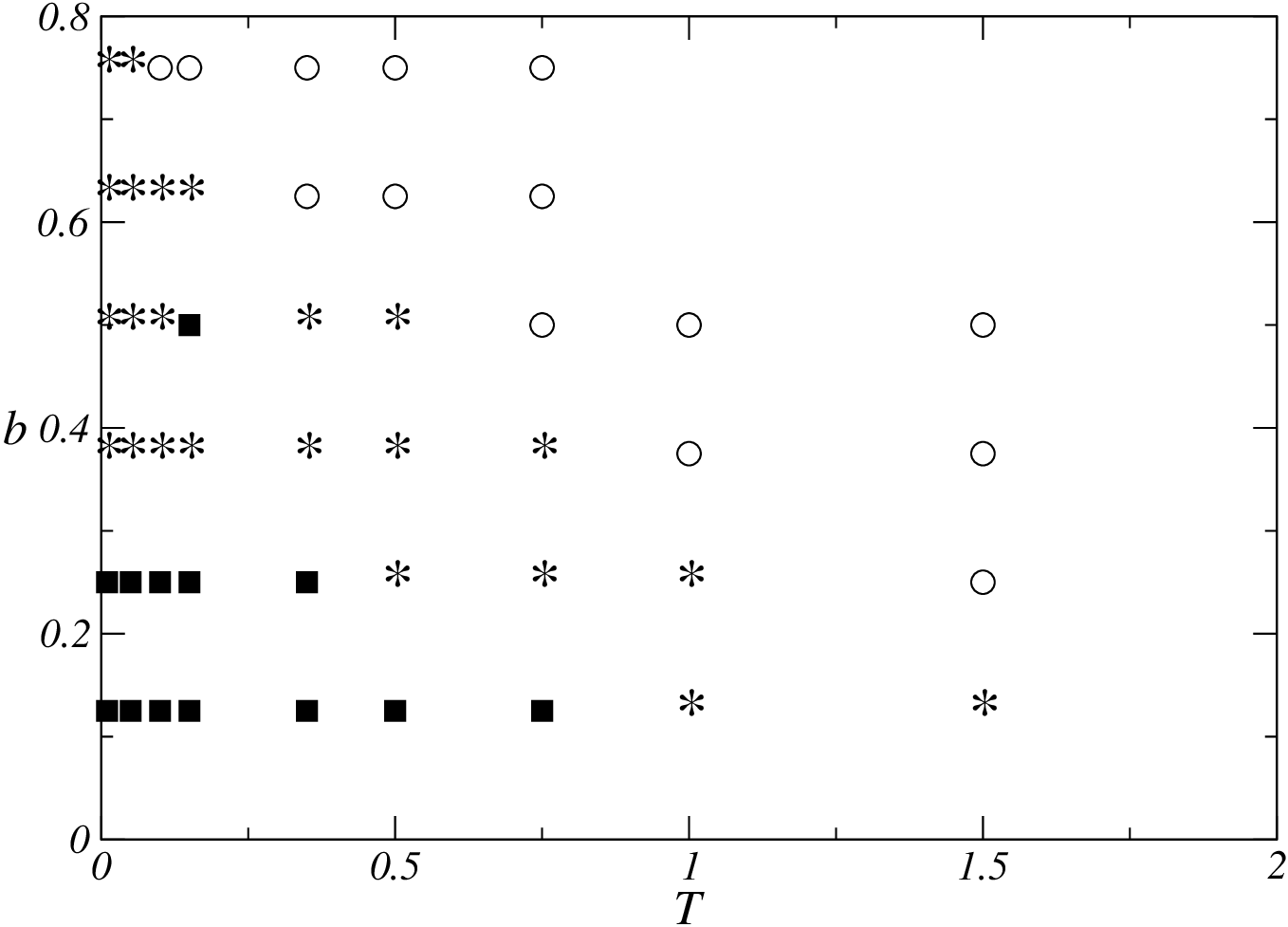}}\,\,\,\,
   \makebox{\includegraphics[scale=0.23]{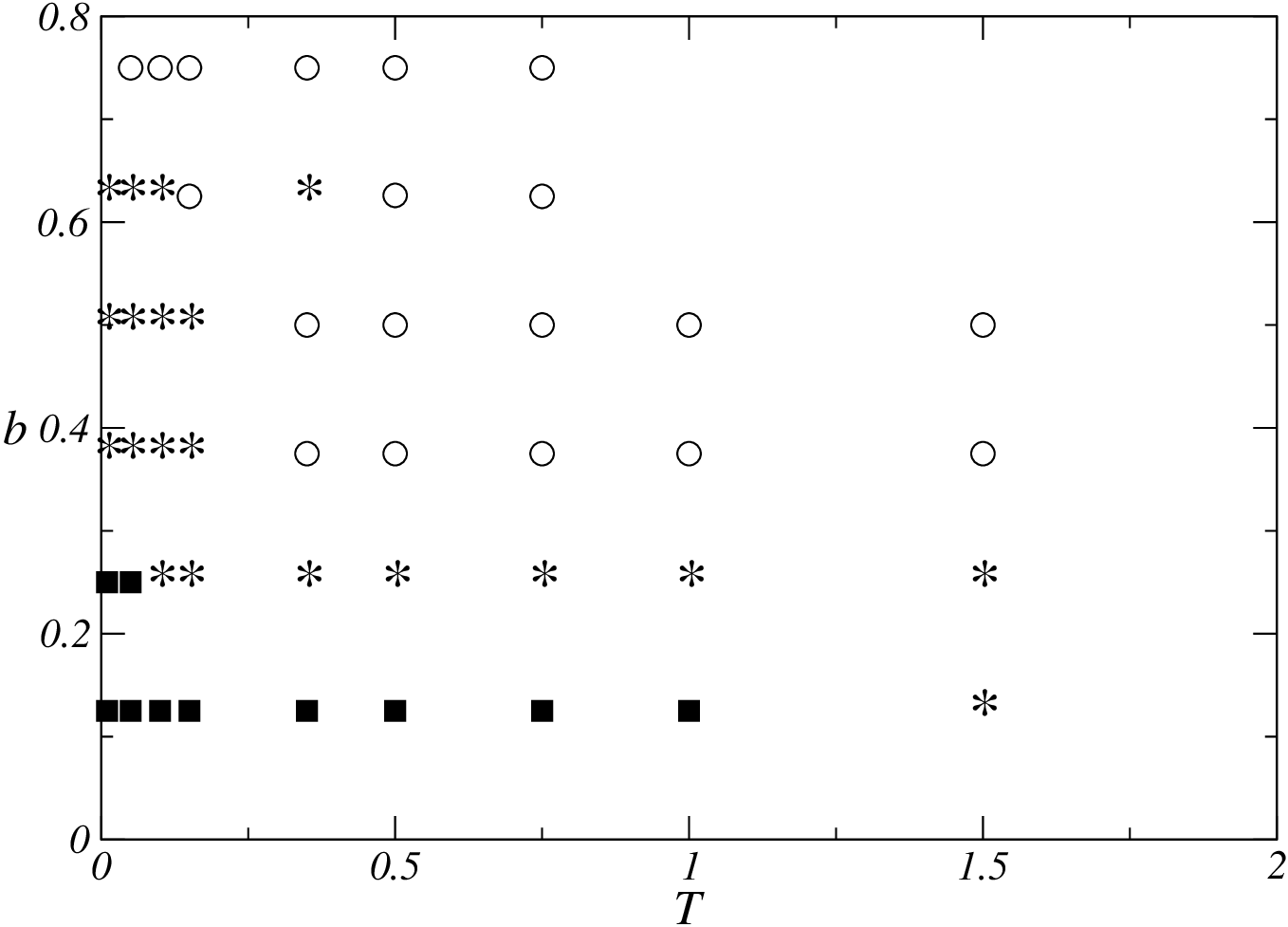}}
  \caption{A part of parameter plane and the corresponding 
$m_{\varepsilon_{thr}}$ function; on the left $m_{12}$, on the right 
$m_{15}$. Squares correspond to Integrable case, stars to Transitional whereas circles to chaotic.}
  \label{fig:effic}
  \end{center}
  \end{figure}
\end{center}

\section{Conclusions}
\label{sec:concl}

In this paper we proved that the simple stirring model given
by~(\ref{bate}) has invariant curves for all $T$ small enough and every
$b\in (0,R)$; the use of NAFF gives us numerical evidence that such invariant
curves persists even for $T$ large, closer and closer to the boundary
of the disk as $T$ increases. Hence, for small $T$ the mixing
efficiency in not good at all, whereas for larger values of $T$, say
$T>1$,  it becomes reasonably good once $b\geq R/2$.
Thus, even if very simple, the model can give rise to efficient
mixing.

We note that, due to the general results on which it relies, our proof of the existence of invariant curves for small $T$  may work for more general $T$--periodic  stirring protocols  $z(t).$
However, the starting point should be to obtain  a 'simple' averaged system.

 From the numerical point of view, the more general systems referred above can be
surely studied by means of NAFF method, once we obtained good
integrators for these systems. Still considering our piecewise
constant model, we think that using NAFF  one could give precise
estimates of the size of the annular domain close to the boundary of $D_R$ containing invariant
curves and of the rate at which it shrinks to zero when
$T$ increases to infinity.

\end{document}